\documentstyle[12pt]{article}

\makeatletter

\newfont{\Bbb}{msbm10 scaled 1\@ptsize00}
\newcommand{\CC}{\mbox{\Bbb C}}

\newif\if@fewtab\@fewtabtrue

\def\draftdate{\number\day.\number\month.\number\year\ \ \ \hourmin }
{\count255=\time\divide\count255 by 60
\xdef\hourmin{\number\count255}
\multiply\count255 by-60\advance\count255 by\time
\xdef\hourmin{\hourmin:\ifnum\count255<10 0\fi\the\count255}}
\def\ps@draft{\let\@mkboth\@gobbletwo
    \def\@oddhead{}
    \def\@oddfoot
       {\hbox to 7 cm{$\scriptstyle\bf Draft\ version:\ \draftdate$
       \hfil}\hskip -7cm\hfil\rm\thepage \hfil}
    \def\@evenhead{}\let\@evenfoot\@oddfoot}

\def\label#1{\ifnum\draftcontrol=1
 \global\def\draftnote{\scriptsize\tt #1}\fi
 \@bsphack\if@filesw {\let\thepage\relax
   \def\protect{\noexpand\noexpand\noexpand}%
\xdef\@gtempa{\write\@auxout{\string
      \newlabel{#1}{{\@currentlabel}{\thepage}}}}}\@gtempa
   \if@nobreak \ifvmode\nobreak\fi\fi\fi
  \@esphack}

\def\@eqnnum{\hbox to 3cm{\phantom{\rm(\theequation)} \draftnote
                         \hfil}\hskip -3cm {\rm(\theequation)}}

\def\eqnarray{\def\draftnote{{}}\global\@fewtabtrue
\stepcounter{equation}\let\@currentlabel=\theequation
\global\@eqnswtrue
\global\@eqcnt\z@\tabskip\@centering\let\\=\@eqncr
$$\halign to \displaywidth\bgroup\@eqnsel\hskip\@centering\@eqcnt\z@
  $\displaystyle\tabskip\z@{##}$&\global\@eqcnt\@ne
  \hskip 1\arraycolsep \hfil${##}$\hfil
  &\global\@eqcnt\tw@ \hskip 1\arraycolsep
$\displaystyle\tabskip\z@{##}$
\hfil  \tabskip\@centering&\global\@eqcnt\thr@@\llap{##}\tabskip\z@
\cr}

\def\endeqnarray{\@@eqncr\egroup
      \global\advance\c@equation\m@ne$$\global\@ignoretrue}

\def\@@eqncr{\let\@tempa\relax
    \ifcase\@eqcnt \def\@tempa{& & &}\or \def\@tempa{& &}
      \or \def\@tempa{&}
      \or\def\@tempa{}
\fi\@tempa
\if@eqnsw
\if@fewtab\@eqnnum\fi
\stepcounter{equation}\fi\global
\@eqnswtrue\global\@eqcnt\z@\global\@fewtabtrue\cr}

\@addtoreset{equation}{section}

\def\cases#1{\left\{\,\vcenter{\normalbaselines\m@th
    \ialign{$\displaystyle{##}\hfil$&\quad##\hfil\crcr#1\crcr}}\right.}

\def\ct#1{\ifnum\draftcontrol=1{\tt [#1]}\else{\cite{#1}}\fi}
\def\ctz#1#2{\ifnum\draftcontrol=1{\tt [#1,#2]}\else{\cite[#1]{#2}}\fi}

\def\draftcite#1{\ifnum\draftcontrol=1#1\else{}\fi}

\def\@lbibitem[#1]#2{\item{}\hskip -3cm \hbox to 2cm
{\hfil$\scriptstyle\draftcite{#2}$}\hskip
1cm[\@biblabel{#1}]\if@filesw
     {\def\protect##1{\string ##1\space}\immediate
      \write\@auxout{\string\bibcite{#2}{#1}}}\fi\ignorespaces}

\def\@bibitem#1{\item\hskip -3cm \hbox to 2cm
{\hfil \scriptsize\tt\draftcite{#1}}\hskip 1cm
\if@filesw \immediate\write\@auxout
       {\string\bibcite{#1}{\the\value{\@listctr}}}\fi\ignorespaces}

\makeatother

\def\llb#1{\label{#1}}
\def\lab#1{\ifnum\draftcontrol=1{{\tt [#1]} \llb{#1}}\else{\llb{#1}}\fi}
\def\Eq#1{(\ref{#1})}
\def\theequation{{\arabic{equation}}}
\def\[{\begin{eqnarray}}
\def\nn{\nonumber}
\def\non{\nonumber \\ }
\def\]{\end{eqnarray}}

\def\een{\end{enumerate}}
\def\ben{\begin{enumerate}}

\renewcommand{\b}{\beta}
\renewcommand{\c}{\gamma}

\renewcommand{\d}{\delta}
\newcommand{\D}{\Delta}
\newcommand{\e}{\epsilon}
\renewcommand{\l}{\lambda}

\def\th{\theta}

\newtheorem{definition}{Definition}

\def\half{\frac{1}{2}}

\def\lie{g}
\def\qlie#1{{\cal L}_q (#1)}  % q-Lie algebra
\def\uq#1{U_q (#1)}
\def\uqg{U_q (\lie)}
\def\ot{\otimes}

\def\wb#1{\wbxxx#1;}  % white Lie bracket
\def\wbxxx#1,#2;{[#1\circ#2]}
\def\bbxxx#1,#2;{[#1\bullet#2]}
\def\t#1{\tilde{#1}}  % tilde
\def\kill#1{\killx#1;}  % Killing form
\def\killx#1,#2;{B(#1,#2)}
\def\killtx#1,#2;{\t{B}\left(#1,#2\right)}

\def\cpppp#1,#2;{#1\cdot#2}
\def\Ch{\CC((\hhh))}
\def\formxx#1,#2;{\langle #1,#2\rangle}
\def\hhh{t}

\newcommand{\rpps}{\left(R_{\pi\pi^*}\right)}
\newcommand{\ripps}{\left(R^{-1}_{\pi\pi^*}\right)}

\newcommand{\cas}{{\cal{C}}}
\newcommand{\tp}{{\cal{T}}}
\newcommand{\tcas}{{\cal B}}
\newcommand{\ttp}{{\cal V}}
\newcommand{\hcas}{K}
\newcommand{\htp}{X}

\newcommand{\gln}{gl_n}
\newcommand{\sln}{sl_n}

\newcommand{\ib}{\bar{\imath}}
\newcommand{\jb}{\bar{\jmath}}
\newcommand{\kb}{\bar{k}}
\newcommand{\lb}{\bar{l}}
\newcommand{\rb}{\bar{r}}
\renewcommand{\sb}{\bar{s}}

\def\remark{\noindent{\bf Remark: }}
\def\paragraph#1{\vspace{3mm}\noindent{\bf #1}}

\begin{document}

\def\draft{\pagestyle{draft}\thispagestyle{draft}
\global\def\draftcontrol{1}}
\global\def\draftcontrol{0}
%%%%%%%%%%%%%%%%%%%% switch on/off draft version %%%%%%%%%%%%%%%%%%%%%%%
%\draft

\newpage
\begin{titlepage}
\begin{flushright}
%{KCL-TH-95-??}\\
YITP/K-1120\\
{q-alg/9508013}
\end{flushright}
\vspace{1.5cm}
\begin{center}
{\bf {\Large QUANTUM LIE ALGEBRAS }\\
\vspace{4mm}
{\large ASSOCIATED TO $U_q(gl_n)$ AND $U_q(sl_n)$}}\\
\vspace{1.2cm}
{\large Gustav W. Delius}
\footnote{On leave from Department of Physics, Bielefeld University,
Germany} and {\large Andreas H\"uffmann}\\
\vspace{3mm}
Department of Mathematics,
King's College London\\
Strand, London WC2R 2LS, Great Britain\\
{\small e-mail: delius@mth.kcl.ac.uk and aha@mth.kcl.ac.uk}\\
\vspace{6mm}
{\large Mark D. Gould }\\
\vspace{3mm}
Department of Mathematics, University of Queensland\\
Brisbane Qld 4072, Australia.\\
\vspace{6mm}
{\large Yao-Zhong Zhang}\\
\vspace{3mm}
Yukawa Institute for Theoretical Physics\\
Kyoto University, Kyoto 606, Japan\\
{\small e-mail: yzzhang@yukawa.kyoto-u.ac.jp}\\
\vspace{1.6cm}
{\bf{ABSTRACT}}
\end{center}
\begin{quote}

Quantum Lie algebras $\qlie{g}$ are non-associative algebras which are
embedded into the quantized enveloping algebras $U_q(g)$ of Drinfeld and
Jimbo in the same way as ordinary Lie algebras are embedded into their
enveloping algebras. The quantum Lie product on $\qlie{g}$ is induced by
the quantum adjoint action of $U_q(g)$.
We construct the quantum Lie algebras associated to $U_q(gl_n)$ and
$U_q(sl_n)$. We determine the structure constants and the quantum
root systems, which are now functions of the quantum parameter $q$.
They exhibit an interesting duality symmetry  under $q\leftrightarrow 1/q$.
\end{quote}
\vfill
\end{titlepage}

\paragraph{1.}
The theory of classically integrable systems relies heavily on Lie algebras
and root systems. The discovery \ct{Dri85,Jim85} of the quantum
deformations
$\uqg$ of the universal enveloping algebras $U(\lie)$ of Lie algebras $\lie$
has led to major advances in the theory of quantum integrable systems.
Many constructions in the theory of classically integrable systems
do however require the use of Lie algebras rather than their
enveloping algebras. An example of this are the values of the
conserved charges on the solitons in affine Toda theory \ct{Fre94}.

To generalize these constructions to the quantum level one would like
to have the concept of a quantum Lie algebra $\qlie{\lie}$ which is related to
the quantized enveloping algebra $\uq{\lie}$ in the same manner as
a Lie algebra is related to its enveloping algebra. Such objects were
introduced in \ct{qlie} and will be reviewed below. In this paper we
construct quantum Lie algebras associated to $\gln$
and $\sln$.

\paragraph{2.}
The quantized enveloping algebra $\uq{\sln}$ is the unital
associative algebra over $\Ch$, the field of fractions for the
ring of formal power series in
the indeterminate $\hhh$, with generators
$x_i^\pm,h_i$, $(i=1,\cdots,n-1)$ and relations
\[\label{uqrel}
&&\left[ h_i,h_j \right] = 0,~~~~
\left[ h_i,x_j^\pm \right] = \pm a_{ij} x_j^\pm ,~~~~
\left[x_i^+ ,x_j^- \right]  =  \delta_{ij}
\ \frac{q^{h_i} - q^{-h_i}}{q- q^{-1}},
\non&&
x^\pm_i x^\pm_i x^\pm_{j}
-(q+q^{-1}) x^\pm_i x^\pm_{j} x^\pm_i
+x^\pm_{j} x^\pm_i x^\pm_i=0~~~~
(|i-j|=1),
\non&&
x^\pm_i x^\pm_j=x^\pm_j x^\pm_i~~~~~~~(|i-j|\geq 2).
\]
where we have defined $q=e^\hhh$. Here
$(a_{ij})$ denotes the Cartan matrix of type $A_{n-1}$, i.e.,
$a_{ii}=2$, $a_{ij}=-1\,(|i-j|=1),=0$(otherwise).
We define $\uq{\gln}$ by adjoining to $\uq{\sln}$ an element
$h_n$ which belongs to the center.
The Hopf algebra structure is given by the coproduct $\D$, the antipode
$S$ and the counit $\e$
\[
&&\Delta(h_i) =h_i \otimes 1 + 1 \otimes h_i, ~~~
\Delta(x_i^\pm) = x_i^\pm \otimes q^{-h_i/2} +
 q^{h_i/2} \otimes x_i^\pm,
\non&&
S(h_i)= - h_i, ~~~
S(x_i^\pm) = - q^{\mp 1}\,x_i^\pm ,~~~
\epsilon(h_i) = \epsilon(x_i^\pm) = 0.
\]
Note that our conventions here differ from those of \ct{qlie} by
$q\leftrightarrow q^{-1}$ in order to conform to \ct{Jim86,Lin92,Gou92}.
The adjoint action of $\uqg$ on itself is given by,
using Sweedler's notation \ct{Swe69},
\[\label{adjoint}
x\circ y=\sum x_{(1)}\,y\,S(x_{(2)}),~~~~~x,y\in\uqg.
\]
The Cartan involution $\theta$ is the algebra automorphism defined by
\[\label{cartinvo}
\th(x_i^\pm)=x_i^\mp,~~~\th(h_i)= -h_i.
\]
It is a coalgebra antiautomorphism, i.e.,
$\D\cdot\th=(\th\ot\th)\cdot\D^T$ and $\e\cdot\th=\th\cdot\e$, and
it satisfies
$S\cdot \th=\th\cdot S^{-1}$.
There is also an involutive algebra antiautomorphism $\dagger :
a\mapsto a^\dagger$ defined by
\[
(x^\pm_i)^\dagger=x^\mp_i,~~~~
(h_i)^\dagger=h_i
\]
which is a coalgebra automorphism and satisfies $S\cdot \dagger=
\dagger\cdot S^{-1}$.
The diagram automorphism $\tau$, defined by
\[
\tau(x^\pm_i)=-x^\pm_{n-i},~~~\tau(h_i)=h_{n-i}~(i\leq n-1),~~~
\tau(h_n)=-h_n,
\]
extends to a Hopf-algebra automorphism.

\paragraph{3.}
A central concept in the theory of quantum Lie algebras \ct{qlie} is
$q$-conjugation which in $\Ch$ maps
$\hhh\mapsto -\hhh$, i.e. $q\mapsto q^{-1}$.
\begin{definition}
a) \underline{q-conjugation}
\mbox{$\sim: \Ch \rightarrow\Ch$}, $a\mapsto\t{a}$ is the field
automorphism defined by $\t{\hhh}=-\hhh$.
\\
\noindent b) Let $M,N$ be $\Ch$-modules. A map
$\phi:M\rightarrow N$ is \underline{$q$-linear} if
$\phi(\l \,a)=\t{\l}\,\phi(a),~\forall a\in M, \l\in\Ch$.
\\
c) Let $A,B$ be algebras over $\Ch$. A q-linear map
$\phi: A\rightarrow B$ is an \underline{algebra $q$-homo-}
\underline{morphism} if it
respects the algebra product, i.e., if
$\forall a,a'\in A,~~\phi(a\,a')=\phi(a)\,\phi(a')$.
$q$-anti-isomorphims, $q$-automorphisms, etc., are defined analogously.
\end{definition}
Note the analogy between the concepts of $q$-conjugation and complex
conjugation and between $q$-linear maps and anti-linear maps.
\begin{definition}\lab{qu}
\underline{$q$-conjugation} on $\uqg$ is the algebra
q-automor\-phism
\mbox{$\sim: \uqg \rightarrow \uqg$} that extends $q$-conjugation on $\Ch$
by acting as the identity on the generators $x_i^\pm$ and $h_i$.
\end{definition}
This definition is consistent because the relations \Eq{uqrel}
are invariant uder $q\mapsto q^{-1}$.
$q$-conjugation is a coalgebra q-antiautomorphism of $\uqg$, i.e.,
$\e \cdot \sim = \sim \cdot \e,~~
\D \cdot \sim = \sim \cdot \D^T$ and it satisfies
$S \cdot \sim = \sim \cdot S^{-1}$.
We define a tilded Cartan involution and a tilded antipode as compositions
\[
\t{S}=\sim\cdot S,~~~~~\t{\th}=\sim\cdot\th,
\]
These behave well with respect to the adjoint action:
\[\label{goodprop}
\t{\th}(a)\circ\t{\th}(b)=\t{\th}(a\circ b),~~~~
\t{S}(a)\circ\t{S}(b)=\t{S}(S^{-1}(a)\circ b),~~~~\forall a,b\in\uqg.
\]

\paragraph{4.}
A Lie algebra $g$ is naturally embedded into its universal enveloping
algebra $U(g)$. It forms a subspace of the enveloping algebra which
under the adjoint action transforms in the adjoint representation.
The Lie bracket on $g$ is given by the restriction of the adjoint action
of $U(g)$.
This idea is extended to the quantum case by the following definition,
introduced in \ct{qlie} (and slightly modified to incorporate also the
non-simple Lie algebra $\gln$):
\begin{definition}\lab{defqlie}
A \underline{quantum Lie algebra} $\qlie{g}$
associated to a finite dimensional
complex Lie algebra $g$ is a $\circ$ - submodule
of $\uqg$, endowed with the
\underline{quantum Lie bracket} $[a\circ\b]=a\circ b$,
and which has the properties:
\begin{enumerate}
\item{$\qlie{g}$ has the same dimension as $g$,}
\item{$\qlie{g}$ is a deformation of $g$, i.e., the isomorphism
$\uqg/\hhh\uqg\cong U(g)$ takes $\qlie{g}/\hhh\qlie{g}$ isomorphically
onto $g\subset U(g)$,}
\item{$\qlie{g}$ is invariant under $\t{\th}$, $\t{S}$ and $\tau$.}
\end{enumerate}
\end{definition}

\remark It was shown in \ct{qlie} that if at least one $\circ$-module
satisfying the first two properties exists, then there exist infinitely
many and out of these
one can always choose at least one statisfying also the last property.
The last property plays a crucial role in the investigations into
the general structure of quantum Lie algebras.

\paragraph{5.}
We will now show how to construct a quantum Lie algebra $\qlie{\gln}$
starting from an expression for the universal R-matrix of $\uq{\gln}$.
Because $\gln$ is not simple the universal R-matrix is not unique and
therefore also the expressions for the quantum Lie algebra generators
which we obtain will not be unique.
Our construction is based on the results of \ct{Lin92,Gou92}.
We introduce the elements $E_{ij}\in\uq{\gln}$
defined recursively by
\[
&&E_{ij}=E_{ik}E_{kj}-q\,E_{kj}E_{ik}~~~
(i<k<j\mbox{ or }i>k>j),~~~~~E_{i,i\pm 1}=x^\pm_i,
\non&&
E_{ii}-E_{i+1,i+1}=h_i~~(i<n),~~~~~
\sum_{i=1}^n E_{ii}=h_n.
\]
%In this paper we use the usual summation convention according
%to which any indices
%which appear in some terms of an expression but not in all terms of
%that expression are summed over their range (here $1\cdots n$).
Let $\pi$ denote the vector representation of $\gln$. Then
$\pi(E_{ij})=e_{ij}~\forall i,j$, where $e_{ij}$ denotes
the matrix $(\d_{ia}\d_{jb})_{1\leq a,b\leq n}$.
We follow the convention of \ct{Lin92,Gou92} and define the dual representation
$\pi^*$ not with the antipode, as would be standard, but
as $\pi^*(a):=\pi^t(\c(a))$ where $\c$ is the
antiautomorphism defined by $\c(x^\pm_i)=-x^\pm_i$, $\c(h_i)=-h_i$.
Rewriting a result by Jimbo \ct{Jim86} one finds that there exists
a universal R-matrix ${\cal R}$ for $\uq{\gln}$ such that
\[
(\pi^*_{ij}\ot 1)({\cal R})=\d_{i\leq j}q^{i-j} \hat{E}_{ij},~~~~~~
(\pi^*_{ji}\ot 1)({\cal R}^T)=\d_{i\leq j}q^{i-j} \hat{E}_{ji},
\]
where we have defined
\[
\hat{E}_{ij}=\left\{
\begin{array}{ll}
-(q-q^{-1})q^{-(E_{ii}+E_{jj}-1)/2} E_{ij} & i\neq j\\
q^{-E_{ii}} & i=j.
\end{array}\right.
\]
(We use a generalized Kronecker delta notation, e.g., $\d_{i\leq j}=1$
if $i\leq j$, $0$ otherwise.)
It is shown in \ct{Lin92,Gou92} that
the elements $T_{ij}\in\uq{\gln}$ defined by
\[
T_{ij}&=&q^{i}(q-q^{-1})^{-1}\left(\d_{ij}-(\pi^*_{ij}\otimes 1)
({\cal R}^T {\cal R})\right)
\non
&=&q^{i}(q-q^{-1})^{-1}\left(\d_{ij}-q^{-i-j}\sum_{k\leq{\rm min}(i,j)}
q^{2k}\hat{E}_{ik}\hat{E}_{kj}\right),
\]
transform as the components of a $\pi\ot\pi^*$ tensor operator, i.e., that
\[\label{Ttransf}
a\circ T_{ij}= T_{kl}\,(\pi_{ki}\ot\pi^*_{lj})\D(a),~~~~
\forall a\in\uq{\gln}.
\]
The representation $\pi\ot\pi^*$ is isomorphic to the adjoint representation
of $\uq{\gln}$.
Classically the $T_{ij}$ go over into the $\gln$ generators $E_{ij}|_{q=1}$.
Thus the $T_{ij}$ span a $\circ$-module satisfying all requirements
of definition \ref{defqlie} for $\qlie{\gln}$ except property 3.

This $\circ$-module also appears as the dual space to the space of
left-invariant
one-forms in the framework of the bicovariant differential calculus on
quantum groups (see e.g.\ \ct{diffcalc}) and it is a special case of the
braided
matrix Lie algebras of Majid~\ct{Majid}.

\paragraph{6.}
The adjoint representation of $gl(n)$ is not irreducible. Correspondingly
the $\circ$-module spanned by the $T_{ij}$ decomposes into a 1-dimensional
module spanned by a Casimir $\cas$ and a $(n^2-1)$-dimensional module
spanned by elements $\tp_{ij}$, where
\[
\cas=\sum_{i=1}^n\frac{1-q^{-2}}{q^{2n}-1}\,q^i\,T_{ii},~~~~~~~~
\tp_{ij}=T_{ij}-\d_{ij}\,q^i\,\cas
\]
%As a basis for it we choose the generators
%$\tp_{i,j}$ with $i\neq j$ and
%\[
%H_i:=\tp_{ii}-q^{-1}\tp_{i+1,i+1}=
%T_{ii}-q^{-1}T_{i+1,i+1},~~~i=1,\dots,n-1.
%\]

\paragraph{7.}
A second $\circ$-module is spanned by the
elements $\tau(T_{ij})$. They transform as follows:
\[
a\circ\tau(T_{ij})&=&\tau(\tau(a)\circ T_{ij})
=\tau(T_{kl})(\pi_{ki}\ot\pi^*_{lj})\D(\tau(a))\nonumber\\
&=&\tau(T_{kl})(\pi^*_{\kb\ib}\ot\pi_{\lb\jb})\D(a).
\]
(Repeated indices are summed over)
We have introduced the notation $\ib=n+1-i$ and have used that
$\D(\tau(a))=(\tau\ot\tau)\D(a)$
and $\pi^*_{ij}(a)=\pi_{\ib\jb}(\tau(a))$.
As a new basis in this second module we choose
\[
{V}_{ij}:=-\tau(T_{\kb\lb})
\rpps_{lkij},
\]
where
\[
\left(R_{\pi\pi^*}\right)_{ijkl}&:=&\left(\pi_{ik}\ot\pi_{jl}^*\right){\cal R}
\nonumber\\
&=&\d_{ik}\d_{jl}+\d_{ij}\d_{kl}
\left((q^{-1}-1)\d_{ik}-(q-q^{-1})q^{k-i}\d_{i>k}\right).
\]
The ${V}_{ij}$ transform in the same way as the $T_{ij}$:
\[\label{tTtransf}
a\circ{V}_{ij}&=&{V}_{kl}\ripps_{klsr}
\left((\pi^*_{rp}\ot\pi_{sq})\D(a)\right)\rpps_{qpij}\\
&=&{V}_{kl}(\pi_{ki}\ot\pi^*_{lj})\D(a),
\]
where we have used the intertwining property of the R-matrix.
They also have the same classical limit.
Again we decompose into the $(n^2-1)$-dimensional and the trivial module
generated by $\ttp_{ij}$ and $\tcas$ respectively, where
\[
\tcas=\sum_{i=1}^n\frac{1-q^{-2}}{q^{2n}-1}\,q^i\,V_{ii}=
-q^{-n}\,\tau(\cas),~~~~~
\ttp_{ij}={V}_{ij}-\d_{ij}q^i\,\tcas
\]

\paragraph{8.}
Because the $T_{ij}$ and the ${V}_{ij}$ have the same
transformation properties, their sums $1/2(T_{ij}+{V}_{ij})$  also generate a
$\circ$-module.
This module contains an $(n^2-1)$-dimensional submodule generated by
${\cal X}_{ij}:=\half\left(\tp_{ij}+\ttp_{ij}\right)$
which is $\tau$-invariant:
\[
\tau({\cal X}_{ij})&=&\half\left(\tau(\tp_{ij})-\tp_{\kb\lb}\rpps_{lkij}
\right)\nonumber\\
&=&-{\cal X}_{\kb\lb}\rpps_{lkij}
-\half\tau(\tp_{\rb\sb})\rpps_{sr\kb\lb}\rpps_{lkij}+\half\tau(\tp_{ij})
\nonumber\\
&=&-{\cal X}_{\kb\lb}\rpps_{lkij}.
\]
For the last equality we used that
$\rpps_{sr\kb\lb}\rpps_{lkij}=
\d_{\rb i}\d_{\sb j}+(q^{-2n}-1)P_{\bar{r}\bar{s}ij}$, where $P$ is
the projector onto the 1-dimensional orbit,
and that $\tp$, which lies in the $(n^2-1)$-dimensional module, vanishes
when contracted with $P$.
As basis vectors we choose
\[
\htp_{ij}&=&
%\mbox{sgn}(i-j)\,(\sqrt{-1})^{i-j}\,
q^{(j-i-1)/2}\,{\cal X}_{ij}~~
(i\neq j=1,\cdots n),
\non
{H}_i&=&{\cal X}_{ii}-q^{-1}{\cal X}_{i+1,i+1}~~
(i=1,\dots,n-1).
\]
%where $\mbox{sgn}(x)$ denotes the sign of $x$.
They satisfy
\[
&&
\t{\th}(X_{ij})=(-1)^{i+j+1}X_{ji},~~~~
\t{\th}(H_i)=-H_i,
\non&&
\t{S}(X_{ij})=-q^{j-i}X_{ij},~~~~~
\t{S}(H_i)=-H_i,
\non&&
(X_{ij})^\dagger=X_{ji},~~~~~(H_i)^\dagger=H_i.
\]

The $\tau$-invariant 1-dimensional module is generated by the Casimir
$\hcas=\cas+q^n\tcas$, which satisfies $\tau(\hcas)=\t{\th}(\hcas)=
\t{S}(\hcas)=-\hcas$, $\hcas^\dagger=\hcas$.
Note that the Casimir $\cas+\tcas$ is not $\tau$ invariant.

Thus we have obtained the quantum Lie algebra $\qlie{\gln}$ spanned by
a basis
$\{X_{ij}|i,j=1,\cdots,n\}\cup\{H_i|i=1,\cdots,n-1\}\cup\{\hcas\}$ which
satisfies all the conditions of definition \ref{defqlie}.

\paragraph{9.}
To calculate the structure constants of $\qlie{\gln}$ we use the formulas
\[
T_{ij}\circ \hat{T}_{kl}&=&(\pi_{rk}\ot\pi^*_{sl})\D(T_{ij})
\non
&=&\frac{q^i}{q-q^{-1}}\left(\d_{ij}\d_{rk}\d_{sl}-
(R_{\pi\pi^*})_{rifa}(R_{\pi^*\pi})_{cfjk}
(R_{\pi^*\pi^*})_{sahb}(R_{\pi^*\pi^*})_{bhcl}\right),
\\
\tau(T_{\jb\ib})\circ \hat{T}_{kl}&=&
(\pi_{rk}\ot\pi^*_{sl})\D(\tau(T_{\jb\ib}))\
=(\pi^*_{\rb\kb}\ot\pi_{\sb\lb})\D(T_{\jb\ib})
\non
&=&\frac{q^{n+1-j}}{q-q^{-1}}\left(\d_{ij}\d_{rk}\d_{sl}-
(R_{\pi^*\pi^*})_{\rb\jb fa}(R_{\pi^*\pi^*})_{cf\ib\kb}
(R_{\pi\pi^*})_{\sb ahb}(R_{\pi^*\pi})_{bhc\lb}\right),
\]
which follow from
$\D(R^T R)=R^T_{12}R^T_{13}R_{13}R_{12}$.
We also need the relation
\[
X_{ii}=-\sum_{k=1}^{n-1}\,q^{i-k}
\frac{q^{2k}-1}{q^{2n}-1}\,{H}_k
+\sum_{k=i}^{n-1}\,q^{i-k}\,{H}_k+q^i\,\hcas.
\]
Tedious but straightforward calculations give
\[\label{stb}
&&\wb{H_k,X_{ij}}=l_{ij}(H_k)\,X_{ij},~~~
\wb{X_{ij},H_k}=-r_{ij}(H_k)\,X_{ij},\non
&&\wb{H_i,H_j}=f_{ij}{}^k\,H_k,~~~~
\wb{X_{ij},X_{ji}}=g_{ij}{}^k\,H_k,\non
&&\wb{X_{ij},X_{kl}}=\d_{jk}\d_{i\neq l}N_{ijl}\,X_{il}+
\d_{il}\d_{j\neq k}M_{kij}\,X_{kj},\\
&&\wb{\hcas,a}=\wb{a,\hcas}=0~~\forall a\in\qlie{\gln}.\nn
\]
where
\[
l_{ij}(H_k)&=&\half\,(1+q^n)
\left(q^{-k}(q\d_{ki}-q^{-1}\d_{k,i-1})
+q^{k-n}(q\d_{k,j-1}-q^{-1}\d_{kj})\right),
\\
r_{ij}(H_k)&=&-l_{ji}(H_k),\label{lr}
\\
f_{ij}{}^k&=&
\d_{ij}\half\left(\d_{k<i}\,(q+q^{-1})(q^k-q^{-k})
+\d_{k>i}\,(q+q^{-1})(q^{n-k}-q^{-n+k})\right.
\nonumber\\
&&\qquad \left.+\d_{ki}(q^{i+1}-q^{-i-1}+q^{n+1-i}-q^{-n-1+i})
\right)
\nonumber\\
&&+\d_{i,j-1}\half\left(\d_{k\leq i}\,(q^{-k}-q^k)
+\d_{k>i}\,(q^{k-n}-q^{-k+n})\right)
\nonumber\\
&&+\d_{i,j+1}\half\left(\d_{k<i}\,(q^{-k}-q^k)
+\d_{k\geq i}\,(q^{k-n}-q^{-k+n})\right),
\\
g_{ij}{}^k&=&\half \,q^{i-j}\left(
\d_{k<j}(q^{k}-q^{-k})
+\d_{k\geq i}(q^{-k}+q^{-k+n})
-\d_{k\geq j}(q^{-k}+q^{k-n})\right),
\\
N_{ijl}&=&\half
%\,\mbox{sgn}(i-j)\mbox{sgn}(j-l)\mbox{sgn}(l-i)
\,q^{-j+1/2}(1+q^n),~~~~
M_{kij}=\t{N}_{kij}\label{ste}
\]
One notices several properties:
\begin{enumerate}
\item{For $t=0~(q=1)$ these are the standard $\gln$ Lie bracket
relations.
}
\item{The ``quantum Cartan subalgebra'' generators $H_i$ have
non-vanishing quantum Lie brackets among themselves. However it
is still commutative in the sense that $\wb{H_i,H_j}=\wb{H_j,H_i}$.
}
\item{There are now two sets of roots, $L=\{l_{ij}\}$ and $R=\{r_{ij}\}$,
related by \Eq{lr}. The combinations $a_{ij}=(l_{ij}+r_{ij})/2$ form
the standard $gl_n$ root lattice, i.e., they satisfy
$a_{ij}+a_{kl}=\d_{jk}a_{il}+\d_{il}a_{kj}$ and
$a_{ji}=-a_{ij}$.
This feature is probably true only for $g$ simply-laced. It is known
that the lattice structure is broken in the non-simply laced cases
of $g=C_2$ \ct{qlie} and $g=G_2$ \ct{unpub}.
}
\item{$\dagger$ is a quantum Lie algebra antiautomorphism, i.e.,
\[
[a^\dagger\circ b^\dagger]=[b\circ a]^\dagger,~~~~~\forall a,b\in\qlie{\lie}.
\]
}
\item{The quantum Lie bracket is $q$-antisymmetric
in the sense that
\[
[a^q\circ b^q]=-[b\circ a]^q,~~~~~\forall a,b\in\qlie{\lie},
\]
where we have defined the $q$-conjugation $a\mapsto a^q$ on
$\qlie{\lie}$ as the $q$-linear map which extends the $q$-conjugation
$\sim$ on $\Ch$ to $\qlie{\lie}$ by acting as the identity on the
basis elements $X_{ij}$, $H_i$ and $K$. (Note that this is not the
same as the q-conjugation $\sim$ on $\uq{\lie}$ which does not leave
$\qlie{\lie}$ invariant.)
This $q$-antisymmetry of the quantum Lie bracket was observed also
for $\qlie{so_5}$ \ct{qlie} and $\qlie{G_2}$ \ct{unpub}.
}
\item{The element $K$ decouples completely. Thus
$\qlie{\gln}$ is not simple.
}
\item{We would like to stress that the relative simplicity of the
structure constants arises only after symmetrizing the $\circ$-module
with respect to the diagram automorphism $\tau$. The quantum adjoint
action of the $T_{ij}$ on themselves is complicated and does not
display any of the features mentioned in the above points 2 to 6.}
\item{The structure constants display the symmetries
$l_{ji}(H_k)=-\widetilde{l_{ij}(H_k)},~$
$r_{ji}(H_k)=-\widetilde{r_{ij}(H_k)}$ and
$f_{ij}{}^k=-\t{f}_{ij}{}^k,~$
which were derived generally in \ct{qlie} from the property \Eq{goodprop}
of $\t{\th}$.}
\item{For $n=3$ these Lie bracket relations reproduce those for
$\qlie{sl_3}$ given in \cite{qlie} (after replacing $q\leftrightarrow 1/q$
and changing the normalization of the generators).}
\end{enumerate}

\paragraph{10.}
To obtain $\qlie{\sln}$ inside $\uq{\sln}$ one can repeat the above
analyis, starting in paragraph 5 with the universal R-matrix of
$\uq{\sln}$ as given by Rosso \ct{Ros89}. Because the numerical
R-matrices (R-matrices evaluated in representations) which we used
for $\gln$ are the same as those for $\sln$, the formulation of the
$\tau$-invariant $(n^2-1)$-dimensional $\circ$-module works as
before and now gives $\qlie{\sln}$. Also the structure constants
for $\qlie{\sln}$ are given by equations \Eq{stb}--\Eq{ste},
simply dropping $K$.

\paragraph{11.}
We define the quantum Killing form $B$ on $\qlie{\sln}$ by
\[
B(a,b) = -q^{-1}\, \mbox{Tr}_\pi\left(\t{S}(a)\,b\,u\right),
\]
where $Tr_\pi$ denotes the trace over the vector representation and
$u$ is the element of $\uq{\sln}$ expressing the square of the antipode
as $S^2(a)=uau^{-1},~\forall a\in\uqg$. This form
is proportional to that defined in \ct{qlie}. It has as its
defining property the ad-invariance
$\kill{a,c\circ b}=\kill{\t{S}(c)\circ a,b}$. It is
$q$-linear in its first argument and
linear in the second and satisfies
$\kill{b,a}=\widetilde{\kill{a,b}}=\kill{\t{S}(a),S(\t{b})}=
\kill{\t{\th}(a),\t{\th}(b)}$. As explained in \ct{qlie}, it is
not the restriction of Rosso's form \ct{Ros90}
on $\uq{\lie}$ to $\qlie{\lie}$. On our basis the Killing form
takes the values
\[&&
\kill{H_i,H_k}=(q+q^{-1})\d_{ik}-\d_{i,k-1}-\d_{i,k+1},
\non&&
\kill{X_{ij},X_{kl}}=\d_{jk}\d_{li},~~~~~~
\kill{H_k,X_{ij}}=0.
\]
The ad-invariance of the Killing form leads to further relations
among the structure constants derived in \ct{qlie}:
\[
g_{ijk}=q^{j-i}r_{ij}(H_k),~~~~
f_{ijk}=f_{ikj},~~~~
N_{kij}=-q^{j-i}N_{ijk}.
\]
Here we have lowered indices with the Killing form, e.g.,
$f_{ijk}=\kill{H_k,H_l}f_{ij}{}^l$.

\paragraph{Acknowledgements}
G.W.D. thanks the Deutsche Forschungsgemeinschaft for a
Habilitationsstipendium.
A.H. thanks the EC for a research fellowship.
Y.Z.Z. is financially supported by the
Kyoto University Foundation.

\end{document}